# Thompson's Group $F$ is Not Minimally Almost Convex

James Belk Kai-Uwe Bux

**Abstract**

We prove that Richard Thompson's group $F$ is not minimally almost convex with respect to the generating set $\{x_0, x_1\}$. This improves upon a recent result of S. Cleary and J. Taback. We make use of the forest diagrams for elements of $F$ introduced by J. Belk and K. Brown. These diagrams seem particularly well-suited for understanding the $\{x_0, x_1\}$-generating set.

Let $G$ be a finitely generated group with a fixed generating set $\Sigma$. For $n \in \mathbb{N}$, let $B_n$ denote the ball of radius $n$ in the Cayley graph of $G$. The convexity function $c : \mathbb{N} \to \mathbb{N}$ of $G$ with respect to $\Sigma$ is defined as follows: Let $d_{B_n}(g, h)$ be the length of a shortest path inside $B_n$ from $g$ to $h$, and put

$$c(n) := \max \left\{ d_{B_n}(g,h) \,\middle|\, g, h \in B_n \text{ and } g^{-1}h \text{ has word length } 2 \right\}.$$

Clearly, $c(n) \geq 2$ for all $n$. On the other hand, any two points in the $n$-ball are connected by a path through the identity, and hence $c(n) \leq 2n$ for all $n$.

**Definition.** We say that $G$ is almost convex with respect to $\Sigma$ if $c(n) \leq C$ for some constant $C$. We say that $G$ is minimally almost convex with respect to $\Sigma$ if $c(n) \leq 2n - 1$ for sufficiently large $n$.

Convexity conditions first arose in the work of J. Cannon, who was interested in algorithms for drawing Cayley graphs of groups. For almost convex groups, he gave an algorithm that constructs the graph $B_n$ for any input $n \in \mathbb{N}$ [Cann87]. Generalizing this result, I. Kapovich [Kapo02] has shown that any minimally almost convex group is finitely presented. T. Riley [Rile02] derives upper bounds for the area function, the isodiametric function, and the filling length of minimally almost convex groups. Upper bounds for the computational complexity of the word problem follow.

Richard Thompson's group $F$ is defined by the presentation

$$F = \left\langle x_0, x_1, x_2, \ldots \,\middle|\, x_{j+1} = x_i^{-1} x_j x_i \text{ for } i < j \right\rangle.$$

Since $x_i = x_0^{1-i} x_1 x_0^{i-1}$ for $i \geq 2$, the group $F$ is generated by the elements $x_0$ and $x_1$.

S. Cleary and J. Taback [ClTa02] have shown that Thompson's group $F$ is not almost convex with respect to $\{x_0, x_1\}$. We prove:

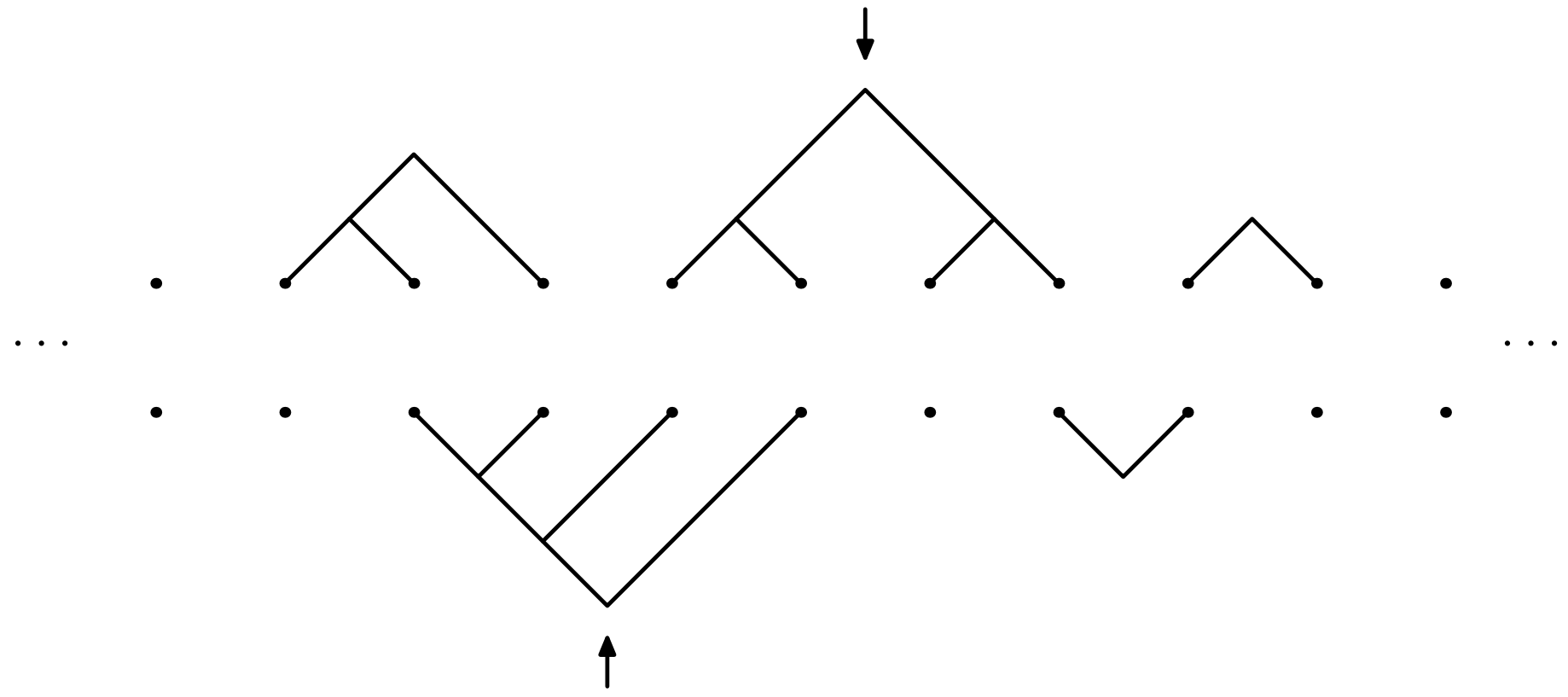

Figure 1: A forest diagram.

**Theorem.** *Thompson's group $F$ is not minimally almost convex with respect to the generating set $\{x_0, x_1\}$. Specifically, for any even $n \geq 4$, there are two elements $f$ and $g$ of length $n$ such that the following conditions hold:*

1. *The elements $f$ and $g$ are distance 2 apart.*
2. *The shortest path connecting $f$ and $g$ inside the $n$-ball has length $2n$.*

# 1 Forest Diagrams

We will be using the forest diagrams for elements of $F$ introduced by J. Belk and K. Brown [BeBr02]. Specifically, we will use the length formula for forest diagrams to calculate distances in the Cayley graph of $F$.

Let $\Gamma$ denote the Cayley graph of $F$. This graph has a vertex for each element of $F$ and an edge from $f$ to $xf$ for every $x \in \{x_0, x_1\}$. The norm $\ell(v)$ of a vertex $v \in \Gamma$ is the distance from $v$ to the identity vertex of $\Gamma$.

Each vertex of $\Gamma$ can be represented by a forest diagram as shown in figure 1. Such a diagram consists of a pair of bounded, bi-infinite binary forests (the top forest and the bottom forest) together with an order-preserving bijection of their leaves.

Let us be a bit more precise about these definitions. A bi-infinite binary forest is a sequence $(\ldots, T_{-1}, T_0, T_1, \ldots)$ of finite binary trees. We can represent such a forest as a line of trees, together with a pointer at the tree $T_0$ (as in figure 1). A forest is bounded if only finitely many of its trees are nontrivial. Note that our binary trees are planar, i.e., the following trees are different:

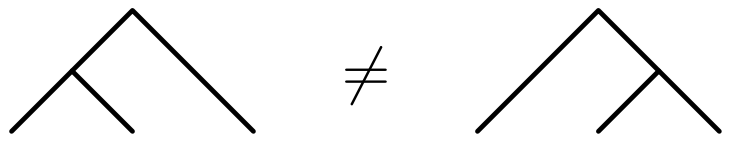

In particular, any binary tree comes with a linear ordering on its leaves; and this in turn induces a linear ordering on the leaves of a bi-infinite binary forest.

A caret is a pair of edges in a forest that join two vertices to a common parent. We call a caret grounded if it joins two leaves. A reduction of a forest diagram consists of removing an opposing pair of grounded carets:

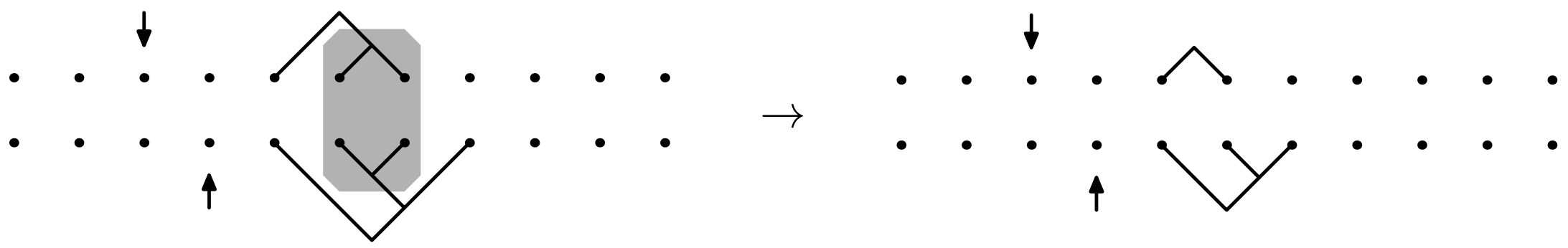

The inverse of a reduction is called an expansion. Two forest diagrams are equivalent if one can be transformed into the other by a sequence of reductions and expansions.

A forest diagram is reduced if it does not have any opposing pairs of grounded carets. It turns out that every forest diagram is equivalent to a unique reduced forest diagram.

**Proposition 1.1 ([BeBr02, Section 4]).** *There is a one-to-one correspondence between vertices of* $\Gamma$ *and equivalence classes of forest diagrams. Therefore, every element of* $F$ *can be represented uniquely by reduced forest diagram.*

**Remark 1.2.** We will frequently identify vertices of $\Gamma$, elements of $F$, and reduced forest diagrams. For example, if $f \in F$, we might talk about the "top forest of $f$". We hope this will not cause any confusion.

The forest diagram for the identity is:

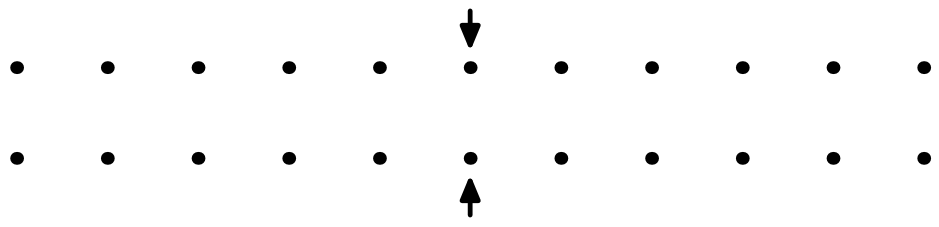

Given a forest diagram for a vertex $v \in \Gamma$, it is easy to find forest diagrams for the neighbors of $v$:

**Proposition 1.3 ([BeBr02, Section 4]).** *Let* $\mathfrak{f}$ *be a reduced forest diagram representing the vertex* $v \in \Gamma$*. Then:*

1. *A forest diagram for* $x_0 v$ *can be obtained by moving the top pointer of* $\mathfrak{f}$ *one tree to the right.*

2. *A forest diagram for* $x_1 v$ *can be obtained by "dropping a caret at the current position". That is, the forest diagram for* $x_1 v$ *can be obtained by attaching a caret to the roots of the top trees in* $\mathfrak{f}$ *indexed by* 0 *and* 1*. Afterward, the top pointer points to the root of the new, combined tree.*

The bottom forest remains unchanged in either case. Note that the given forest diagram for $x_1 v$ will need to be reduced if the new caret opposes a grounded caret from the bottom tree. In this case, left-multiplication by $x_1$ effectively deletes a grounded caret from the bottom tree.

**Corollary 1.4.** *Again, let $\mathfrak{f}$ be a reduced forest diagram for a vertex $v \in \Gamma$. Then:*

1. *A forest diagram for $x_0^{-1}v$ can be obtained by moving the top pointer of $\mathfrak{f}$ one tree to the left.*

2. *A forest diagram for $x_1^{-1}v$ can be obtained by deleting the top caret of the current tree. The top pointer ends at the resulting left-child tree. If the current tree is trivial, one must first perform an expansion. In this case, left-multiplication by $x_1^{-1}$ effectively creates a new grounded caret in the bottom tree.*

**Example 1.5.** Here are some sample edges from the Cayley graph of $F$:

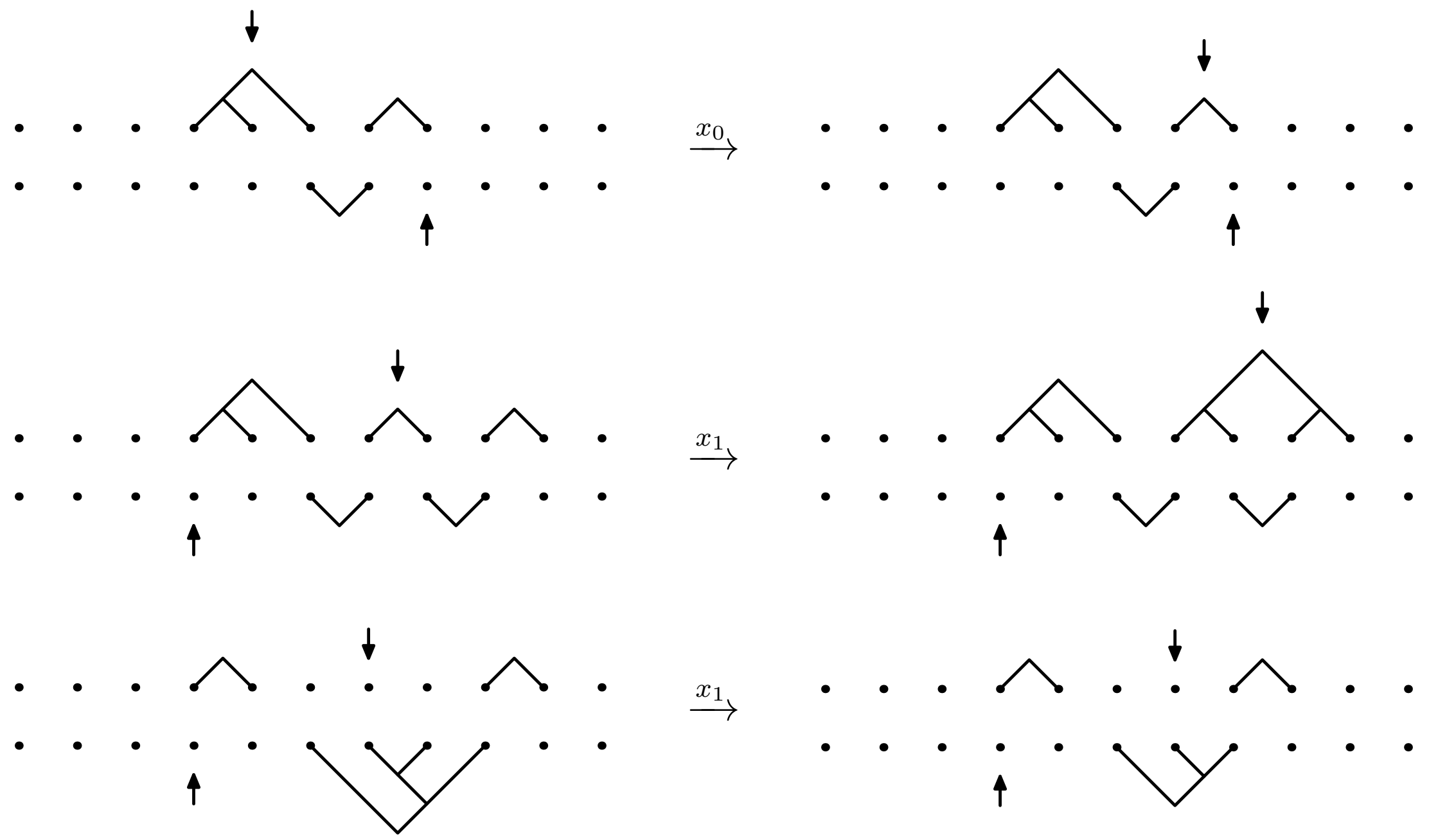

An element of $F$ is called positive if it lies in the submonoid generated by $x_0$ and $x_1$. An element $f \in F$ is semi-positive if its reduced forest diagram has a trivial bottom forest. Equivalently, $f$ is semi-positive if it lies in the submonoid generated by $\{x_0, x_0^{-1}, x_1\}$.

**Observation 1.6.** *Let $f$ be a semi-positive element with forest diagram $\mathfrak{f}$, and let $v$ be a vertex of $\Gamma$ with forest diagram $\mathfrak{g}$. Then a forest diagram for $fv$ can be obtained by stacking $\mathfrak{f}$ on top of $\mathfrak{g}$. That is, the top forest of $\mathfrak{f}$ is attached along its leaves to the roots of the top forest of $\mathfrak{g}$ in such a way that the bottom pointer of $\mathfrak{f}$ and top pointer of $\mathfrak{g}$ match up.*

**Example 1.7.** Put:

$f :=$ 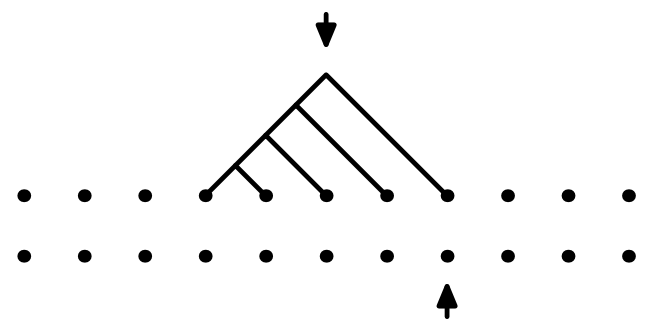 and $v :=$ 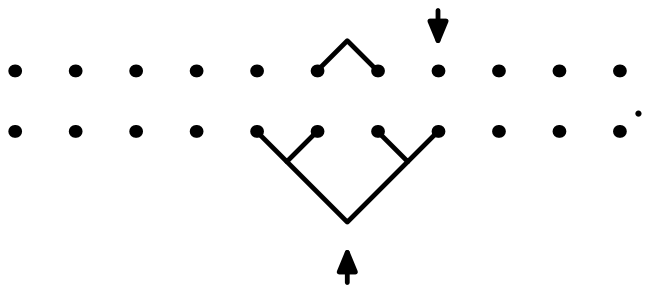.

Then:

$fv =$ 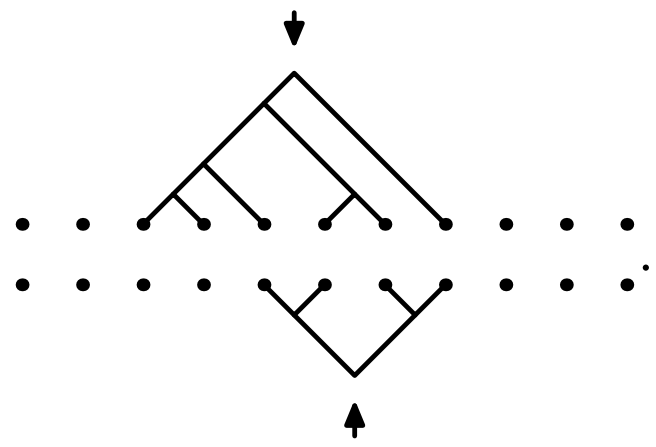.

Note that this diagram is reduced: none of the bottom carets of $v$ ended up opposing the unique grounded caret of $f$.

**Remark 1.8.** Forest diagrams relate to a certain action of $F$ on $\mathbb{R}$ by piecewise-linear homeomorphism in the same way that the tree diagrams of [CFP96] relate to the standard action of $F$ on $[0, 1]$. In particular, $F$ can be regarded as the group of piecewise-linear homeomorphisms $f$ of $\mathbb{R}$ satisfying the following conditions:

1. Each linear segment of $f$ has slope an integral power of 2.

2. The graph of $f$ has only finitely many breakpoints, each of which has dyadic rational coordinates.

3. The function $f$ is eventually of the form
$$f(t) = t+K_+$$
as $t \to \infty$ and:
$$f(t) = t+K_-$$
as $t \to -\infty$ for some integers $K_+, K_-$.

Given such a homeomorphism $f \in F$, there exist subdivisions of its domain and range into standard dyadic intervals such that $f$ maps intervals of the domain linearly to intervals of the range. The bottom forest of a forest diagram corresponds to this partition of the domain, and the top forest corresponds to the partition of the range. Details on this can be found in [BeBr02, Section 4].

## 2 Computing the Length of an Element

Since the action of $x_0$ and $x_1$ is relatively simple, it comes as no surprise that one can find the length $\ell(f)$ of an element $f \in F$ directly from a forest diagram. Our treatment of the length formula is based on [BeBr02].

We begin with some terminology. A space is the region between two leaves in a forest. A space is interior if it lies between two leaves from the same tree, and exterior if it lies between two trees. Note that every exterior space in a forest is either to the left or the right of the pointer.

Given a forest diagram for an element $f \in F$, we label the spaces between the leaves of each forest as follows. Label a space:

- **L** (for *left*) if it exterior and to the left of the corresponding pointer,
- **N** (for *necessary*) if it is not of type **L** and if the leaf to the right of the space is a left leaf in its caret,
- **I** (for *interior*) if it is interior and not of type **N**, or
- **R** (for *right*) if it is exterior, to the right of the corresponding pointer, and not of type **N**.

See figure 2 for an example.

The spaces of a forest diagram come in pairs: one from the top forest and one from the bottom forest. The support of a forest diagram is the minimum interval that contains both pointers and all nontrivial trees. We only label space pairs in the support of a forest diagram.

The weight of a space pair is determined by the following table:

| | **N** | **I** | **R** | **L** |
|---|---|---|---|---|
| **N** | 2 | 2 | 2 | 1 |
| **I** | 2 | 0 | 0 | 1 |
| **R** | 2 | 0 | 2 | 1 |
| **L** | 1 | 1 | 1 | 2 |

We can now state the length formula for elements of $F$. We follow the exposition in [BeBr02, Section 5], which is a simplification of Fordham's original formula [Ford95]. Viewing Thompson's group as a diagram group, V. Guba [Guba02] has recently obtained a different version of the length formula.

**Theorem 2.1 (Length Formula).** *Let $f \in F$, and let $\mathfrak{f}$ be its reduced forest diagram. Let $\ell_1(f)$ be the total number of carets of $\mathfrak{f}$, and let $\ell_0(f)$ be the sum of the weights of all space pairs in the support of $\mathfrak{f}$. Then $f$ has length*

$$\ell(f) = \ell_1(f) + \ell_0(f)\,.$$

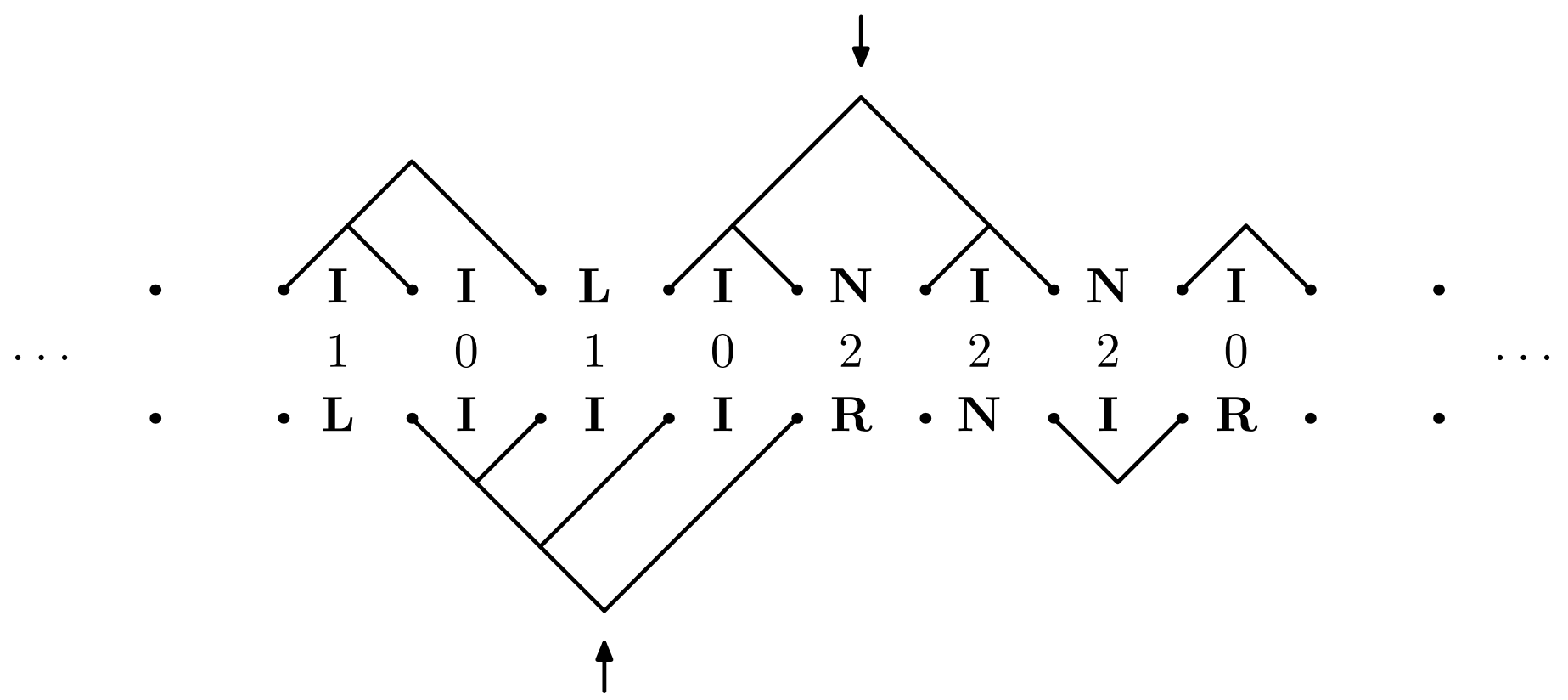

Figure 2: An element of length $18 = 6 + 4 + 8$. A minimum-length word is $x_0^{-1}x_1x_0x_1x_0^{-1}x_1x_0^{-1}x_1^{-1}x_0^2x_1x_0x_1^2x_0^{-1}x_1^{-3}$.

**Remark 2.2.** In fact, there exists a minimum-length word $w$ for $f$ with the following properties:

1. The number of $x_1$'s in $w$ is equal to the number of carets in the top forest of $\mathfrak{f}$.
2. The number of $x_1^{-1}$'s in $w$ is equal to the number of carets in the bottom forest of $\mathfrak{f}$.
3. The total number of $x_0$'s and $x_0^{-1}$'s in $w$ is equal to $\ell_0(f)$.

However, not every minimum-length word for $f$ has this form. In particular, it is sometimes possible to exchange an $\{x_0, x_0^{-1}\}$-pair for an $\{x_1, x_1^{-1}\}$-pair.

Despite this complication, we call $\ell_1(f)$ the $\underline{x_1\text{-count}}$ of $f$, and $\ell_0(f)$ the $\underline{x_0\text{-count}}$ of $f$.

**Example 2.3.** Figure 2 shows an element $f \in F$ with its space pairs labeled. The weight of each space pair is also indicated, as is a minimum-length word for $f$. Note that the minimum-length word has six occurrences of $x_1$, four occurrences of $x_1^{-1}$, and eight total occurrences of $x_0$ and $x_0^{-1}$.

# 3 The Main Theorem

**Theorem 3.1.** *For any $n \geq 1$, there exist elements $l, r \in F$ such that:*

1. *The distance from $l$ to $r$ in the Cayley graph is* 2.
2. *Both $l$ and $r$ have length $2n + 2$.*
3. *Any path from $l$ to $r$ inside the $(2n + 2)$-ball has length at least $4n + 4$.*

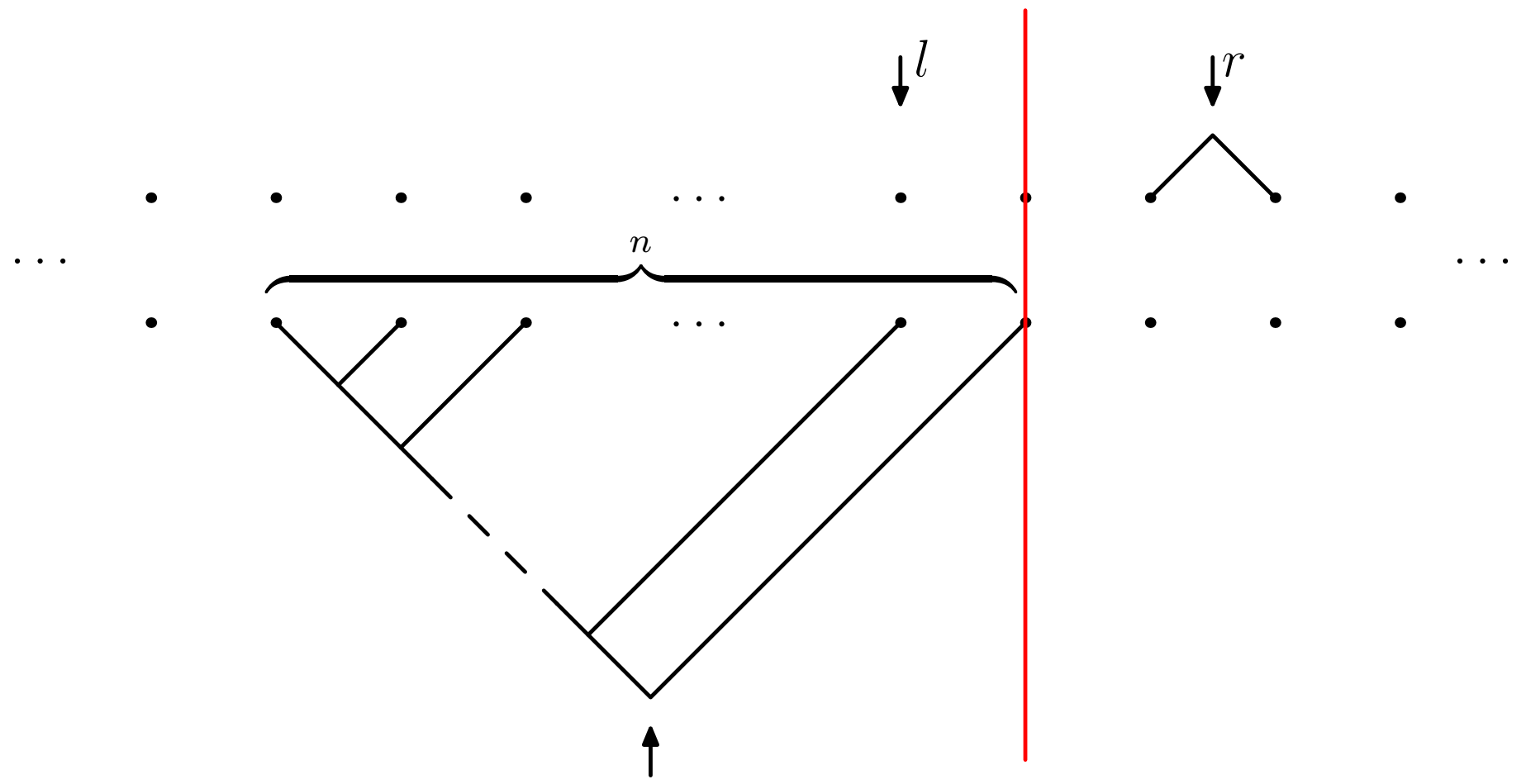

Figure 3: The forest diagrams for $l$ and $r$. These two elements differ only in the position of the top pointer, as indicated. The vertical line is “critical”.

Here are the elements in question:

$$\begin{aligned} l &:= x_0^{-2} x_1 x_0^{n+1} x_1^{-n} \\ r &:= x_0^2 l \end{aligned}$$

The forest diagrams for $l$ and $r$ have the same carets, and differ only in the position of the top pointer. We have shown both forest diagrams in figure 3, with the two possible positions of the top pointer indicated.

Clearly $l$ and $r$ are distance 2 apart, with geodesic path

$$l \to x_0 l \to r.$$

**Lemma 3.2.** *The elements $l$ and $r$ both have length $2n+2$.*

**Proof.** We begin with $l$. Note first that the forest diagram for $l$ has exactly $n+1$ carets. Furthermore, the space pairs of $l$ have the following labels:

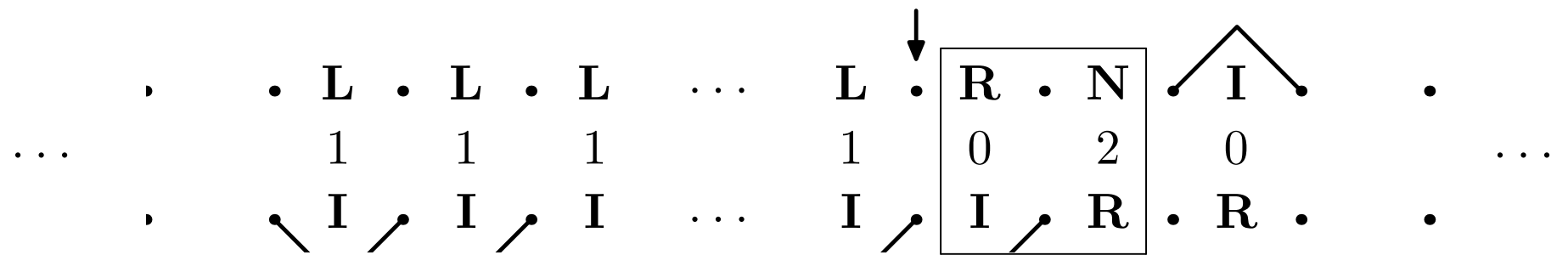

Thus $l$ has length $(n+1) + \underbrace{1 + \cdots + 1}_{n-1} + 0 + 2 + 0 = 2n+2$.

Similarly, $r$ has exactly $n+1$ carets. Its labels are:

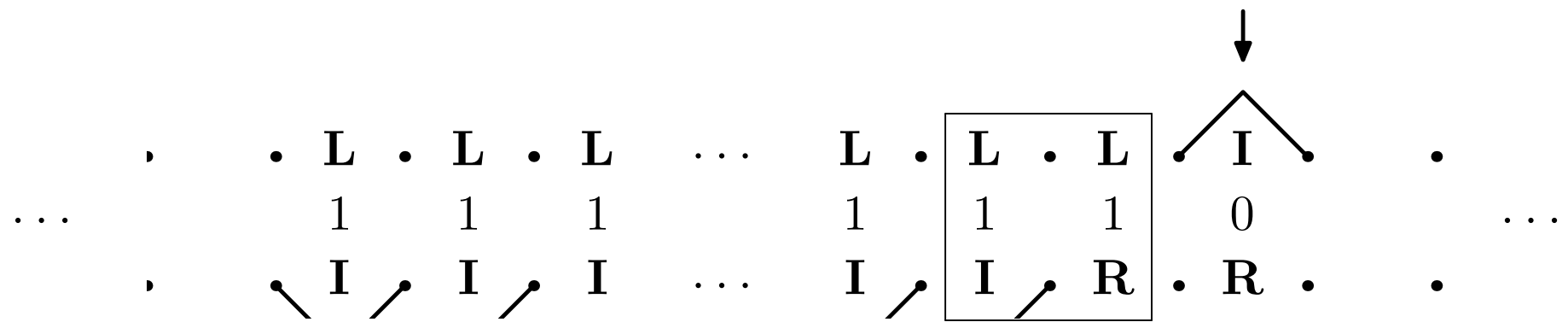

Therefore, $r$ has length $(n+1) + \underbrace{1 + \cdots + 1}_{n-1} + 1 + 1 + 0 = 2n + 2$. **q.e.d.**

**Remark 3.3.** Note that the element $x_0 l$ has labels

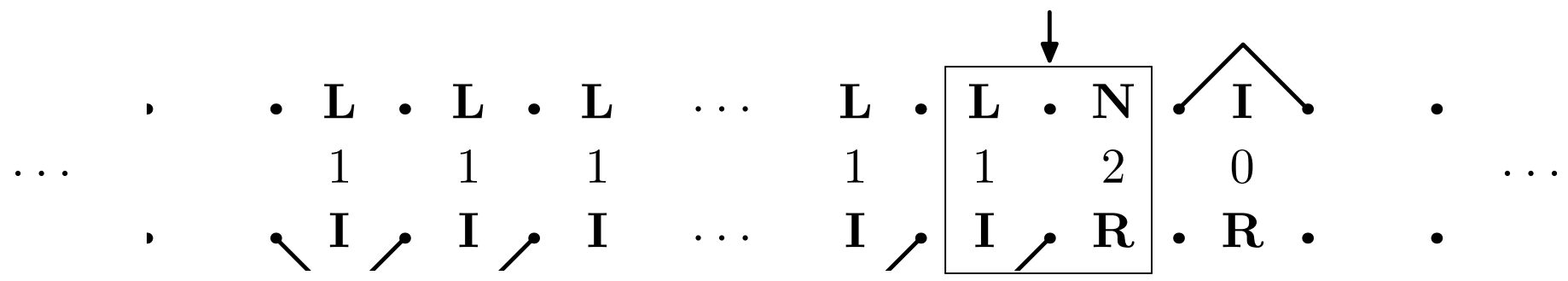

and hence has length $2n+3$. Therefore, the geodesic path

$$l \to x_0 l \to r$$

leaves the ball of radius $2n+2$.

We postpone the proof of condition (3) to the next section. In the remainder of this section, we discuss the intuitive ideas behind the proof and try to convey how $l$ and $r$ were chosen.

The main idea of the proof is as follows. Consider the "critical line" in the forest diagram for $l$ and $r$ (the vertical line in figure 3). This line has the following crucial property: if one remains in the $(2n+2)$-ball, it is not possible for the top pointer to cross the critical line while both of the outermost carets are in place. Therefore, any path in the $(2n+2)$-ball from $l$ to $r$ must go through the following four stages:

1. Move to the left, and delete the leftmost caret.

2. Move to the right (crossing the critical line), and delete the rightmost caret.

3. Move back left (crossing again), and re-create the leftmost caret.

4. Move back right (crossing the critical line for a third time), and re-create the rightmost caret.

**Example 3.4.** The word

$$\left(x_1 x_0^{n+1}\right)\left(x_1^{-1} x_0^{-n}\right)\left(x_1^{-1} x_0^{n}\right)\left(x_1 x_0^{1-n}\right).$$

describes a path in $B_{2n+2}$ from $l$ to $r$ of length $4n+4$. (Recall that the edges in the Cayley graph $\Gamma$ correspond to *left*-multiplication by generators. Hence, this word should be read from *right to left* when interpreted as a set of instructions for movement in $\Gamma$.)

Note that we left the bulk of the bottom tree intact throughout this path. In particular, this path does not pass through the identity vertex.

**Example 3.5.** The word

$$\left(x_1 x_0^{n+1}\right)\left(x_1^{-n} x_0^{-1}\right)\left(x_1^{-1} x_0 x_1^{n-1}\right)\left(x_1 x_0^{1-n}\right)$$

represents a minimum-length path from $l$ to $r$ that passes through the identity vertex. This time, we "travel to the right" by destroying the bottom tree, and "travel to the left" be re-creating it.

**Example 3.6.** For $n = 8$, here is another minimum-length path from $l$ to $r$:

$$\left(x_1 x_0^{7} \boldsymbol{x_1^{-4}} x_0^{2}\right)\left(x_1^{-1} x_0^{-4}\right)\left(x_1^{-1} x_0^{4}\right)\left(x_1 x_0^{-2} \boldsymbol{x_1^{4}} x_0^{-5}\right).$$

In this path, we build carets in the top forest while moving to the left, and destroy them later during the final move to the right. (Note that we have highlighted the segments of the word under discussion.) The resulting "bridge" saves us travel time during the middle two stages, but its construction and demolition cost the same amount of time during the first and last stages.

Finally we give an example of two vertices that do *not* work. This example shows why we have chosen elements that have the bulk of their support to the left of the critical line: There is a subtle difference in moving to the right and moving to the left.

**Example 3.7.** Consider the elements

$$\begin{aligned} l' &:= x_1^{-1} x_0^{-n} x_1 x_0^{n-1} \\ r' &:= x_0^{2} l' \end{aligned}$$

Forest diagrams for $l'$ and $r'$ are pictured in figure 4. Clearly $l'$ and $r'$ have distance 2 in $\Gamma$, and it is easy to check that they both have length $2n+2$. Furthermore, the path

$$l' \rightarrow x_0 l' \rightarrow r'$$

leaves the $(2n+2)$-ball. This suggests a "critical line" in the forest diagram (already shown in figure 4), and led us to believe that

$$\left(x_0^{1-n} x_1\right)\left(x_0^{n+1} x_1^{-1}\right)\left(x_0^{-n} x_1^{-1}\right)\left(x_0^{n} x_1\right) \qquad \text{(length } 4n+4\text{)}$$

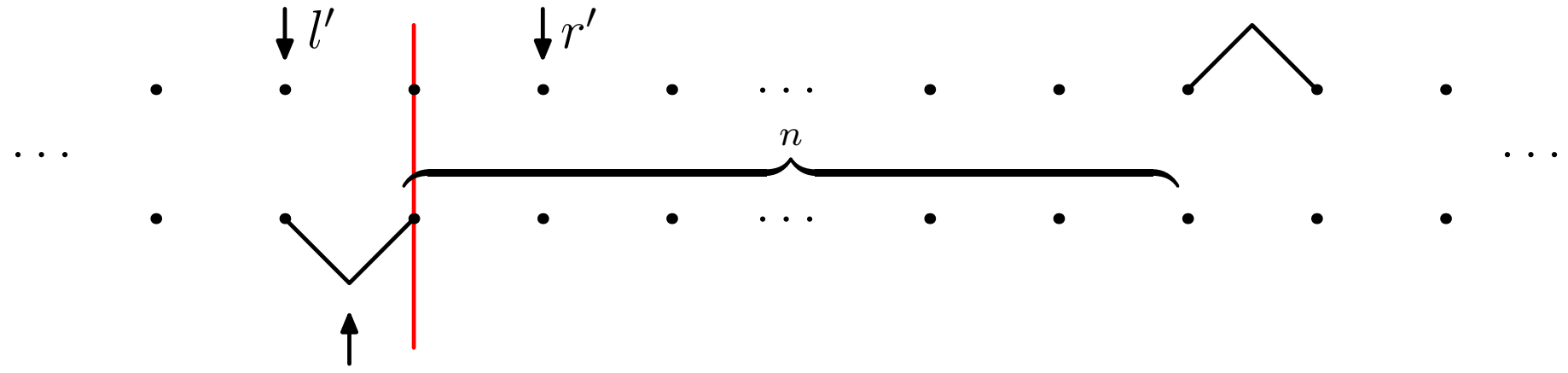

Figure 4: Two elements that fail. The vertical line pretends to be critical.

is a minimum-length path from $l'$ to $r'$ in $B_{2n+2}$.

However, this turns out not to be the case. For example, when $n = 8$,

$$\left(\boldsymbol{x_1^{-6}} x_0^{-1} x_1\right) \left(x_0^3 x_1^{-1}\right) \left(x_0^{-2} x_1^{-1}\right) \left(x_0^n \boldsymbol{x_1^6} x_0 x_1\right)$$

is a path from $l'$ to $r'$ of length 24 in $B_{18}$. We save time in this path by building a "bridge" during our first move to the right. This does not cost any extra time, since building the bridge helps us to get to the right. We then use this bridge during the next two stages to move back an forth very quickly. Finally, we tear the bridge down during our final move to the left, which again does not cost any extra time.

# 4 Proof of Condition (3.1.3)

Fix a path $p$ from $l$ to $r$ that does not leave $B_{2n+2}$. We wish to show that its length $L(p)$ is at least $4n + 4$.

We claim that it suffices to show:

**Lemma 4.1.** *On the path $p$, there are two vertices, $h_l$ and $h_r$, such that*

$$d_\Gamma(h_l, h_r) \geq 2n + 3.$$

Why is this sufficient? Well clearly,

$$L(p) \geq d(l, h_l) + d(h_l, h_r) + d(h_r, r)\,.$$

However, by the triangle inequality

$$d(h_l, h_r) \leq d(h_l, l) + d(l, r) + d(r, h_r)\,.$$

Since $d(l, r) = 2$, we conclude that

$$L(p) \geq 2\,d(h_l, h_r) - 2 \geq 4n + 4.$$

It remains to prove lemma 4.1. We begin by formalizing the notion of "crossing the critical line" from the previous section. If $v$ is a vertex in $\Gamma$, the two pointers of $v$ point to two trees in the forest diagram: one top tree and one bottom tree. We call the rightmost leaf of the top tree the right foot of $v$, and the rightmost leaf of the bottom tree the critical leaf of $v$. Note that the right foot of $l$ is to the left of the critical leaf and that the right foot of $r$ is to the right of the critical leaf.

Let $h_l$ be the first vertex of $p$ whose right foot steps on the critical leaf, and let $h_r$ be the last vertex of $p$ with this property.

**Remark 4.2.** Along the path $p$, there must exist a vertex whose right foot lies directly above its critical leaf. Note, however, that this statement is not entirely trivial. In particular, moving along an edge of $\Gamma$ can change the position of the right foot by more than one unit.

**Lemma 4.3.** *The path $p$ ends with*

$$x_0^{-1}x_1^{-1}r \rightarrow x_1^{-1}r \rightarrow r.$$

*In particular, $h_r = x_0^{-1}x_1^{-1}r$.*

**Proof.** It is easy to check that every path of length three emanating from $r$ either passes through $x_0^{-1}x_1^{-1}r$ or leaves the $(2n+2)$-ball. **q.e.d.**

An element $f \in F$ is left-sided if (1) the pointers in $f$ point to trees whose rightmost leaves match, and (2) all trees to the right of the pointers are trivial.

The width $w(f)$ of an element $f \in F$ is the number of space pairs in the support of the reduced forest diagram for $f$.

**Lemma 4.4.** *If $f \in F$ is left-sided, then*

$$\ell(f) \geq 2\,w(f)\,.$$

**Proof.** We can associate to each caret in a forest diagram the interior space that it covers. Therefore, the interior spaces of $f$ each contribute 1 to the $x_1$-count of $f$. However, since $f$ is left-sided, every exterior space of $f$ is of type **L**. The claim now follows, since the weight of any space pair is greater than or equal to the number of **L**'s in its pair of labels. **q.e.d.**

**Remark 4.5.** Note that lemma 4.4 *fails* for right-sided elements: An $\{\mathbf{R}, \mathbf{I}\}$-space pair does not contribute to the $x_0$-count, and only contributes 1 to the $x_1$-count. Hence, the best available estimate for right-sided elements is $\ell(f) \geq w(f)$.

This difference is related to the fact that one can "move right" by dropping carets, but one cannot simultaneously build a structure and move left. (Compare with example 3.7.)

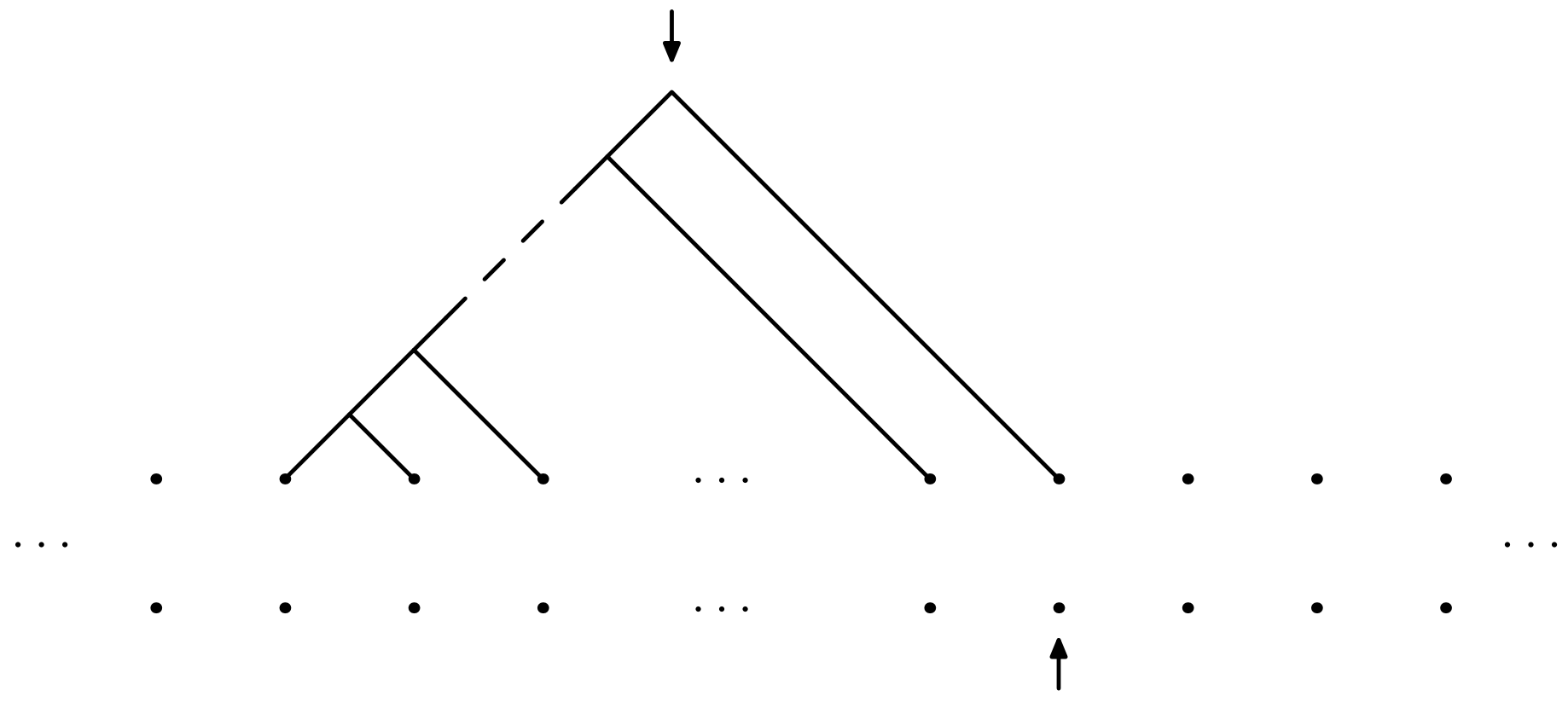

Figure 5: The forest diagram for $h_r^{-1}$.

**Proof of Lemma 4.1.** Note that $h_r^{-1}$ is semi-positive (see figure 5). Therefore, for any $v \in \Gamma$, a forest diagram for $h_r^{-1}v$ is obtained by stacking $h_r^{-1}$ on top of $v$. Moreover, this diagram will be reduced unless the bottom forest of $v$ has a grounded caret in exactly the right position (namely, $n$ spaces to the left of the critical leaf) to oppose the unique grounded caret of $h_r^{-1}$.

Put $x_2 := x_0^{-1}x_1x_0$. Note:

- Every left-sided element commutes with $x_2$. In particular, $h_r^{-1}$ and $x_2$ commute.
- We have $\ell(x_2 f) = \ell(f) + 3$ for any left-sided $f \in F$.

Now, $h_l$ is the first vertex of $p$ whose right foot hits the critical leaf. Therefore, when we get to $h_l$ in $p$, we have not yet modified any material to the right of the critical leaf. In particular, there is some left-sided $h_l' \in F$ satisfying

$$h_l = x_2 h_l'.$$

Observe that $\ell(h_l') = \ell(h_l) - 3 \leq 2n - 1$, and hence $h_l'$ has width strictly less than $n$. Then the stacked diagram for $h_r^{-1}h_l'$ must already be reduced, since no caret of $h_l'$ is far enough to the left to oppose the grounded caret of $h_r^{-1}$. From this, we conclude that $h_r^{-1}h_l'$ has width at least $n$. Since $h_r^{-1}h_l'$ is left-sided, lemma 4.4 implies:

$$\ell\bigl(h_r^{-1}h_l'\bigr) \geq 2n$$

and hence

$$d(h_l, h_r) = \ell\bigl(h_r^{-1}h_l\bigr) = \ell\bigl(x_2 h_r^{-1}h_l'\bigr) = \ell\bigl(h_r^{-1}h_l'\bigr) + 3 \geq 2n + 3. \qquad \textbf{q.e.d.}$$